\documentclass[11pt]{article}

\def\a{\alpha}
\def\b{\beta}
\def\d{\delta}

\def\s{\sigma}

\def\T{{\cal T}}
\def\Tau{{\bf{\cal T}}}
\def\ot{\otimes}

\def\lra{\longrightarrow}

\def\mapright#1{\smash{\mathop{\lra}\limits^{#1}}}

\def\RM{{\mathbf R}}
\def\NM{{\mathbf N}}

\def\ZM{{\mathbf Z}}

\newtheorem{thm}{Theorem}
\newtheorem{de}[thm]{Definition}
\newtheorem{lem}[thm]{Lemma}
\newtheorem{prop}[thm]{Proposition}
\newtheorem{cor}[thm]{Corollary}

\begin{document}

\title{\bf Twisted descent algebras and the Solomon-Tits algebra}

\author{{\bf Fr\'ed\'eric Patras}\\
CNRS, UMR 6621\\
Parc Valrose\\
06108 Nice cedex 2\\
France\\
patras@math.unice.fr\\
 \\
{\bf Manfred Schocker\thanks{supported by Deutsche Forschungsgemeinschaft (DFG-Scho 799)}
}\\
Mathematical Institute\\
24-29 St Giles'\\
Oxford OX1 3LB\\
United Kingdom\\
schocker@maths.ox.ac.uk\\}

\date{}
\maketitle

\section*{Introduction}
The purpose of the present article is to define and study a new class
of descent algebras, called twisted descent algebras.  These algebras
are associated to the Barratt-Joyal theory of twisted bialgebras in
the same way than classical descent algebras are associated to
classical bialgebras.  The theory of twisted descent algebras is a
refinement of the theory of descent algebras of twisted algebras, as
introduced in \cite{pr3}. The formal properties of twisted descent
algebras seem particularly meaningful in view of applications to
discrete probabilities, to the geometry of Coxeter groups and
buildings, and to symmetric group combinatorics.\par Let us survey
briefly the classical theory, since most of the results that hold in
that case have natural generalizations in the new setting. Recall
first that classical shuffles and descent classes in the symmetric
groups are one of the building blocks of modern algebraic
combinatorics.  Reutenauer's monograph \cite{r} contains an overview
of the development until 1993, and strongly indicates the impact of
the so-called Solomon or descent algebra on the theory.  This
classical descent algebra, introduced by Solomon in 1976 in the more
general setting of finite Coxeter groups, was originally constructed
as a noncommutative superstructure of the character ring of the
underlying group, and built upon the classical character theory
encapsulated in the celebrated Mackey formulae \cite{s}.\par Shuffles
and descents are in fact central in various fields of mathematics,
such as discrete probabilities (card shuffling, their asymptotic
properties, and so on), algebraic topology (where the shuffles were
soon recognized as the key ingredient for constructing products on the
singular cohomology of spaces; other operations such as Steenrod
operations or operations on other cohomology theories such as
equivariant or generalized cohomologies can often be described using
symmetric group combinatorics and generalizations of shuffles), and
the theory of iterated integrals and special functions (see
\cite{do,pa,ch,p2,ra} for some examples).\par In the current point of
view, descent algebras are part of a research field including the
theory of noncommutative symmetric functions, quasi-symmetric
functions, and their various generalizations \cite{mr,g,k}. Moreover,
a descent algebra can be associated to any commutative or
cocommutative bialgebra \cite{pt,p}.  In this general setting,
Solomon's classical descent algebra is recovered in the particular
case when the underlying bialgebra is the shuffle algebra over a
countable alphabet.\par \medskip Almost simultaneously to Solomon's
discovery, in 1977, Barratt introduced a new class of algebras,
twisted algebras, in order to study a subtle class of invariants in
the cohomology of topological spaces known as (generalized) Hopf
invariants \cite{ba}.  The theory of twisted algebras was developed
further by Joyal, who showed that the usual classes of algebras (Lie
algebras, enveloping algebras, bialgebras, and so on) can be defined
consistently in this new framework \cite{joy,st}.  However, the
importance of twisted bialgebras for algebraic combinatorics was
uncovered only recently. One of the most striking applications of the
general theory so far is probably the construction of a natural
enveloping algebra structure on the direct sum of the symmetric group
algebras, deduced from the properties of the free twisted algebra on
one generator \cite{pr3}.\par \medskip To combine the two notions of
twisted algebras, and descent algebras, was the motivation for the
present article.  In the first two sections, it is shown that a
twisted descent algebra can be associated naturally to any twisted
cocommutative or commutative bialgebra --- that is, to a cocommutative
or commutative bialgebra in the monoidal symmetric category of
(linear) tensor species. In Section~3, we show that each of the new
twisted objects carries a second associative product, the composition
product.  The free twisted descent algebra is defined by generators
and relations in Section~4, and maps onto any twisted descent algebra.
More precisely, this free object behaves with respect to twisted
descent algebras of twisted bialgebras exactly as the ordinary descent
algebra behaves with respect to descent algebras of bialgebras.
For example, it carries two products and a coproduct, just like the
ordinary one.\par Various combinatorial identities in twisted descent
algebras are established along the way. This includes in particular a
natural generalization (and explanation) of the crucial multiplicative
reciprocity law for descent algebras, which (in the classical case)
was originally discovered in \cite{g} and links both products and, in
the free case, the coproduct on these algebras (see
Proposition~\ref{remarkable} and Corollary~\ref{reciproci}).
Besides, it turns out that the theory of twisted descent algebras is
related to the Solomon-Tits algebra in the same way than the descent
algebras of bialgebras are related to the classical descent algebra.\par
Recall that the Solomon-Tits algebra was introduced by Tits in an
appendix to Solomon's original paper \cite{ti}. As a vector space over
the ground field, it is generated by the elements of the Coxeter
complex of the symmetric group $S_n$, where $n$ is a fixed integer.
The product of two elements of the Coxeter complex is defined in terms
of the minimal galleries joining them (in Tits' buildings
terminology).  The construction makes sense for general hyperplane
arrangements and, in this generality, has been one of the leading
tools in recent research on random walks associated to such
arrangements and related topics, such as the Tsetlin library, or the
random-to-top process in computer science \cite{bid,br}.\par We feel
that the heavy reference to geometry, together with the sometimes
concise style of the appendix, was one of the reasons why Tits' ideas
concerning descents, unlike Solomon's, have not yet been fully
exploited by combinatorists. In that sense, our constructions and
proofs unravel the combinatorial ideas hidden in Tits' line of
reasoning.  This is further demonstrated in final Section~5.  However,
the objects defined and studied here are much richer than the original
one: Tits associated to a Coxeter complex an associative algebra
structure; our twisted descent algebras (the Solomon-Tits algebra of
the symmetric group being a particular case of which) carry two
different associative algebra structures, and, in the free case, a
bialgebra structure. Besides, there is a fundamental rule linking all
three structures.

\section{Tensor species and twisted bialgebras}
 
Recall that a (linear) tensor species is a functor from the category
of finite sets $\bf Fin$ (and set isomorphisms) to the category $\bf
Vect$ of vector spaces over a field $\bf k$ (and linear isomorphisms)
or, more generally, to the category of modules over an arbitrary
commutative ring.  Concretely, a tensor species $F$ associates to each
finite set $S$ a vector space $F(S)$. A bijection from $S$ to $T$ in
$\bf Fin$ induces an isomorphism from $F(S)$ to $F(T)$. We will assume
from now on that the finite sets we consider (and therefore the
objects of $\bf Fin$) are subsets of a given countable ordered set,
for example the set $\NM$ of positive integers. This is not a serious
restriction, but a convenient one for notational purposes. We also
assume that our tensor species are connected, that is: $F(\emptyset
)=\bf k$.

\begin{de} 
  A tensor species is a functor $F$ from $\bf Fin$ to $\bf Vect$ such
  that $F(\emptyset )=\bf k$, where $\bf k$ is the ground field.
\end{de}

\begin{prop}
  The category $\bf Sp$ of tensor species is a linear symmetric
  monoidal category for the tensor product defined by:
  $$
  (F\ot G)( S):=\bigoplus_{T\coprod U=S}F({T})\ot G( {U})
  $$
  for all $F,G\in{\bf Sp}$.
\end{prop} 

Here, $\coprod$ stands for disjoint union: that is, we have $T\cap
U=\emptyset$ in the above formula.\par The proof is straightforward
and is left to the reader. The unit of the tensor product is the
ground field $\bf k$, identified with the tensor species (also denoted
by $\bf k$) with unique nontrivial component $\bf k(\emptyset):=\bf
k$:
$$
F\ot{\bf k}={\bf k}\ot F=F
$$
for all $F\in{\bf Sp}$.

\begin{de}
  A twisted algebra is an algebra in the linear symmetric monoi\-dal
  category of tensor species.
\end{de}

See \cite{joy,pr3} for further details on and references to Joyal's
theory of twisted algebras. Concretely, a twisted algebra is a tensor
species $F$ provided with a product map (which is a map of tensor
species: $F\ot F\lra F$). Associative algebras, commutative algebras,
Lie algebras and so on, are defined accordingly.\par For example, a
twisted associative algebra with unit is a tensor species $A$ provided
with a product map
$$A\ot A\mapright{m}A$$
such that associativity holds: $$m\circ (m\ot
A)=m\circ (A\ot m).$$
Here, we write $A$ for the identity morphism of
$A$. The unit condition is defined in the same way: $A$ has to be
provided with a unit map ${\bf k}\to A$ satisfying the
usual identities. The fundamental example of a twisted associative
algebra is the free twisted algebra ${\bf S}$ on one generator: if
${\bf k}[1]$ denotes the tensor species defined by ${\bf k}[1](\{
n\}):={\bf k}$ for all $n\in\NM$, and ${\bf k}[1](S):=0$ whenever $S$
is not a singleton, then ${\bf S}=\oplus_{n\in\NM}{{\bf k}[1]}^{\ot
  n}$. The product map is the obvious one:
$$
m:{{\bf k}[1]}^{\ot n}\ot {{\bf k}[1]}^{\ot m}\mapright{~}{{\bf
    k}[1]}^{\ot n+m}.
$$
The following examples illustrate the meaning of this construction
and show how to handle computations with twisted algebras in
practice. Let us start with the two-element set ${S}=\{ 3,5\}$. Then,
by definition:
\begin{eqnarray*}
{\bf S}( {S})
& = &
({{\bf k}[1]}\ot {{\bf k}[1]})( \{ 3,5\} )\\[1mm]
& = &
\Big({{\bf k}[1]}(\{3\})\ot {{\bf k}[1]}(\{5\})\Big)
\oplus
\Big({{\bf k}[1]}(\{5\})\ot {{\bf k}[1]}(\{3\})\Big)
  = 
{\bf k}\oplus{\bf k},
\end{eqnarray*}
that is, ${\bf S}( {S})$ is isomorphic to the direct sum of two copies
of the ground field indexed by the two sequences $(3,5)$ and $(5,3)$.
If $ {U}=\{2,4\}$ and ${\phi}$ is the unique order preserving map
from $ {S}$ to $ {U}$, then ${\bf S}( {\phi} )$ maps the component of
${\bf S}( {S})$ indexed by $(5,3)$ to the component of ${\bf S}( {U})$
indexed by $(4,2)$, and so on. In such a situation, we will say that
${\bf S}( {\phi} )$ is induced by a map on the indices.\par More
generally, given any set $ {S}=\{s_1,\ldots,s_n\}$ of order $n$, ${\bf
  S}( {S})$ is isomorphic to the direct sum of $n!$ copies of the
ground field, canonically indexed by all the permutations of the
sequence $(s_1,\ldots,s_n)$, that is by all the sequences $(s_{\s
  (1)},\ldots,s_{\s (n)})$, where $\s\in S_n$, the symmetric group on
$\{1,\ldots,n\}$. The product map $m$ from ${\bf S}( {S})\ot {\bf S}(
{T})$ to ${\bf S}( {S}\coprod {T} )$ is induced by the map on the
indices:
$$
((s_{\s (1)},\ldots,s_{\s (n)}),(t_{\tau (1)},\ldots,t_{\tau(k)}))
\longmapsto 
(s_{\s (1)},\ldots,s_{\s (n)},t_{\tau(1)},\ldots,t_{\tau(k)}),
$$
where $\s$ and $\tau$ run over all permutations in $S_n$,
respectively $S_k$. Notice that the product map in the twisted algebra
${\bf S}$ therefore bears an apparent resemblance to the usual
concatenation product of words. The fine algebraic structure of the
free twisted algebra on one generator, however, is quite different 
from the algebraic structure of the algebra of words in an
alphabet (the free associative algebra over this alphabet).\par\ \par
The notion of a twisted coalgebra is dual to the notion of a twisted
algebra: a twisted coassociative coalgebra with counit is a tensor
species $C$ provided with a coproduct map
$$
C\mapright{\d}C\ot C
$$
such that coassociativity holds: 
$$
(\d\ot C)\circ \d=(C\ot \d )\circ\d ,
$$
as well as the usual identities for the counit.  Once again, the basic
example is provided by the twisted algebra ${\bf S}$, which has a
natural coalgebra structure:
$$
\d :{\bf S} ( {S})
    \longrightarrow 
    \bigoplus\limits_{ {U}\coprod {T}={S}}{\bf S}( {U})\ot {\bf S}( {T}).
$$
(see \cite{pr3} for another approach and details on the structure of
${\bf S}$, and for generalities on twisted algebras or coalgebras.)
As above, it is enough to specify the action of $\d$ at the level of
the indices: whenever $U\coprod T=S$, the $(U,T)$-component of
$\d({\bf S} ( {S}))$ is obtained by sending isomorphically the $(s_{\s
  (1)},\ldots,s_{\s (n)})$-component of ${\bf S} ( {S})$ to the component
of ${\bf S}(U)\ot {\bf S}(T)$ indexed by $u=(u_{\a(1)},\ldots,u_{\a
  (p)})$ and $t=(t_{\tau (1)},\ldots,t_{\tau (q)})$, where $u$
(respectively, $t$) is the unique ordered subsequence of $(s_{\s
  (1)},\ldots,s_{\s (n)})$ composed of the elements of $ {U}$
(respectively, of $ {T}$). For example, if $ {S}=\{ 2,3,5\}$, the
coproduct $\d$ on the $(5,2,3)$-component of ${\bf S} ( {S})$ is
induced by the following map of indices:
\begin{eqnarray*}
  (5,2,3)
  & \longmapsto &
  \{
  ((5,2,3),\emptyset ),\; ((5,2),(3)),\;((5,3),(2)),\;((2,3),(5)),\\[1mm]
  &&
  \hspace*{1ex}
  ((5),(2,3)),\;((2),(5,3)),\;((3),(5,2)),\;(\emptyset ,(5,2,3))\}.
\end{eqnarray*}
The coproduct $\d$ is clearly cocommutative. It is related to the
unshuffle coproduct on the tensor algebra in the same way than the
product is related to the concatenation product. It is shown in
\cite{pr3} that $\d$ induces a (new) coproduct on the direct sum of the
symmetric group algebras, turning this direct sum into an enveloping
algebra.\par\ \par The tensor product (in the category of tensor
species) of two twisted associative algebras $A$ and $A'$ is an
associative algebra. The product map is defined, as in the usual case,
using the symmetry isomorphism:
$$
A\ot A'\cong A'\ot A,
$$
that reads pointwise:
\begin{eqnarray*}
  (A\ot A')[ {S}]
  & = &
  \bigoplus\limits_{ {T}\coprod {U}= {S}} A[{T}]\ot A'[ {U}]\\[1mm]
  & \cong &
  \bigoplus\limits_{ {U}\coprod {T}= {S}} A'[{U}]\ot A[{T}]
  =
  (A'\ot A)[ {S}],
\end{eqnarray*}
where the isomorphism in the middle is the usual symmetry isomorphism
for the tensor product of vector spaces.

\begin{de} 
  A bialgebra in the category of tensor species, or twisted bialgebra,
  is a tensor species $B$ carrying a unital associative twisted
  algebra structure (with product written $m$) together with a
  coassociative counital coalgebra structure (with coproduct written
  $\d$) such that $\d :B\to B\ot B$ is a map of twisted algebras.
\end{de}
 
The tensor species $\bf S$ is provided with the structure of a twisted
bialgebra, by the product $m$ and the coproduct $\d$ defined above.

\section{Twisted descent algebras}

Recall that, given a set $I$, an $I$-graded vector space $V$ is a
collection of vector spaces indexed by the elements of $I$:
$$
V=\{ V_i\}_{i\in I},\quad V_i\in {\bf Vect} .
$$
It is often convenient to identify $V$ with the direct sum of its
graded components : $V=\bigoplus_{i\in I}V_i$. A morphism $\phi$ of
graded vector spaces from $V$ to $W$ is a family of morphisms:
$$
\phi_i :V_i\lra W_i
$$
or, equivalently, a morphism from $V=\bigoplus_{i\in I}V_i$ to
$W=\bigoplus_{i\in I}W_i$ such that $\phi=\bigoplus_{i\in I}\phi_i$
with $\phi_i\in {\rm Hom}(V_i,W_i)$.  We denote by ${\rm Hom}_I(V,W)$
the set of all these morphisms, and set ${\rm End}_I(V):={\rm
  Hom}_I(V,V)$.\par 
To each tensor species $F$ is canonically associated a $\cal P$-graded
vector space also written $F$, where ${\cal P}$ denotes the set of
finite subsets of $\NM$.  For each $X\in \cal P$, we set $F_X:=F(X)$.
The category ${\cal P}$-${\bf Vect}$ of ${\cal P}$-graded vector
spaces carries two tensor products $\overline \ot$ and $\ot$: the
homogeneous one, defined by
$$
(V\overline\ot W)_X=V_X\ot W_X
$$
for all $X\in {\cal P}$,
and the linear species or ''$\cal P$-graded'' one, defined by
$$
(V\ot W)_X=\bigoplus\limits_{T\coprod V=X}V_T\ot W_U
$$
for all $X\in {\cal P}$, as in Section~1.

\begin{de} 
  Let $B$ be a twisted bialgebra. The convolution product on the
  vector space of $\cal P$-graded linear endomorphisms of $B$ is
  defined by
  $$
  f\ast g:\ B\mapright{\d}B\ot B\mapright{f\ot g}B\ot B\mapright{m}B.
  $$
  for all $f,g\in {\rm End}_{\cal P}(B)$.
\end{de}
Here, the tensor product $f\ot g:B\ot B\to B\ot B$ is defined by:
$$
(B\ot B)_S
=
\bigoplus\limits_{T\coprod V=S}B_T\ot B_V
\unitlength1ex
\begin{picture}(24,6)
\put(2,0.60){\line(1,0){16}}
\put(17.75,0){\mbox{$\longrightarrow$}}
\put(3.50,4.00){\mbox{$\bigoplus\limits_{T\coprod V=S}f_T\ot g_V$}}
\end{picture}
=
(B\ot B)_S.
$$
It is important to notice that we \emph{do not} assume that $f$
and $g$ are endomorphisms of $B$ as a tensor species. Therefore,
strictly speaking, the tensor product of $f$ and $g$ is not induced by
the symmetric monoidal structure of $\bf Sp$ (but by the one in 
${\cal P}$-${\bf Vect}$).
 
\begin{prop} 
  The convolution product is associative and unital. In particular,
  $({\rm End}_{\cal P}(B),\ast )$ is an associative algebra with
  unit.
\end{prop} 

Indeed, for $f,g,h\in {\rm End}_{\cal P}(B)$, we have:
\begin{eqnarray*}
  f\ast (g\ast h)
  & = &
  m\circ (f\ot (g\ast h))\circ\d\\[1mm]
  & = &
  m\circ (f\ot (m\circ (g\ot h)\circ \d ))\circ\d\\[1mm]
  & = &
  m\circ (B\ot m)\circ (f\ot g\ot h)\circ (B\ot \d )\circ \d
  \mbox{\hphantom{$=(f\ast g)\ast h$}}
\end{eqnarray*}
and, since $m$ (respectively, $\d$) is associative (respectively,
coassociative):
\begin{eqnarray*}
  \mbox{\hphantom{$f\ast (g\ast h)$}}
  & = &
  m\circ (m\ot B)\circ (f\ot g\ot h)\circ (\d \ot B)\circ \d 
  =
  (f\ast g)\ast h.
\end{eqnarray*}

\begin{de} 
  The characteristic map $1_S$ associated to $S\in\cal P$ is defined
  by:
  \begin{eqnarray*}
    & 1_S : \ B_{S\,} \mapright{=} B_{S\,}, & \mbox{and}\\[2mm]
    & 1_S : \ B_{S'}  \mapright{0} B_{S'}   & \mbox{if } S\neq S'.    
  \end{eqnarray*}
\end{de}
The unit $\eta$ of the convolution is the characteristic map
$1_\emptyset$ associated to $\emptyset$, as can be checked easily.

\begin{de} 
  The twisted descent algebra ${\cal T}_B$ of a twisted bialgebra $B$
  is the convolution algebra generated by the characteristic maps of
  the finite subsets of $\NM$.
\end{de}

The twisted descent algebra is generated, as a vector space, by
$1_\emptyset$ and by the convolution products $1_{S_1}\ast \cdots\ast
1_{S_k}$, where $S_i\cap S_j=\emptyset\neq S_i$ for all $i$, $j$ such
that $i\not= j$.  For, if $T\cap U\not=\emptyset$, it follows from the
definitions that $1_T\ast 1_U=0$.  We write from now on
$1_{(S_1,\ldots,S_k)}$ for $1_{S_1}\ast \cdots\ast 1_{S_k}$. In
general, the elements $1_{(S_1,\ldots,S_k)}$, even with the
restrictions on the sets $S_i$ given above, are not linearly
independent.\par\ \par Since $m$ is associative, it induces a unique
product map $m^{[3]}$ from $B^{\ot 3}$ to $B$, given by:
$m^{[3]}:=m\circ (m\ot B)=m\circ (B\ot m)$. More generally, we write
$m^{[k]}$ for the unique product map from $B^{\ot k}$ to $B$ and
(since the coproduct is coassociative) $\d^{[k]}$ for the unique
coproduct map from $B$ to $B^{\ot k}$.  Then, in particular, we have:
$$
1_{(S_1,\ldots,S_k)}
= 
m^{[k]}\circ 1^\ot_{(S_1,\ldots,S_k)}\circ \d^{[k]},
$$
where
$1^\ot_{(S_1,\ldots,S_k)}:=1_{S_1}\ot \cdots\ot 1_{S_k}$.\par\ \par 
Notice that the (ordinary) descent algebra of $B$ has already been
defined in \cite{pr3}: it is, in the language of the present article,
the convolution algebra generated by the formal series
$1_n:=\sum\limits_{S\in {\cal P}, \ |S|=n}1_S$ ($n\in\NM$).

\section{The fine algebraic structure of twisted descents}

The twisted descent algebra of a twisted algebra $B$ carries, by its
very definition, an associative product, the convolution product. In
the present section, we show that it is also closed under the
composition product in the category of $\cal
P$-graded vector spaces, a property that extends to the twisted
setting one of the fundamental properties of the Solomon descent
algebra and, more generally, of descent algebras of commutative or
cocommutative bialgebras. Indeed, they also carry two products, that,
together with certain coalgebra properties, are the building blocks of
the modern theory of descent algebras. We restrict our attention to
the case where $B$ is cocommutative. Dual properties hold when $B$ is
commutative.\par

\begin{thm} \label{main}
  The twisted descent algebra of a twisted cocommutative bialgebra $B$
  is closed under the composition of $\cal P$-linear endomorphisms of
  $B$. The composition product is given by:
  \begin{eqnarray*}
    \lefteqn{%
      1_{(  S_1,\ldots,  S_n)} \circ 1_{(  T_1,\ldots,  T_k)}}\\[2mm]
  & = &
  \left\{
    \begin{array}{ll}
      0 & \mbox{if }
          S_1\coprod\ldots\coprod S_n\neq T_1\coprod\ldots\coprod T_k,\\[2mm]
      1_{
        (S_1\cap T_1,\ldots,S_1\cap T_k,
         \;\ldots\;,
         S_n\cap T_1,\ldots,S_n\cap T_k),
       }&
       \mbox{otherwise.}
     \end{array}\right.
  \end{eqnarray*}
\end{thm}

Note that, in the second part of the formula, possibly occuring empty
intersections $S_i\cap T_j$ may be omitted, since $1_\emptyset$ is the
identity of ${\cal T}_B$. For example, if
$(S_1,S_2)=(\{3,5\},\{1,4\})$ and $(T_1,T_2)=(\{5\},\{1,3,4\})$, then
$$
1_{(S_1,S_2)}\circ 1_{(T_1,T_2)}
=
1_{(\{5\},\{3\},\emptyset,\{1,4\})}
=
1_{(\{5\},\{3\},\{1,4\})}.
$$\par
The two identities in the theorem are the most fundamental identities
satisfied in twisted descent algebras. Since the first identity
follows immediately from the definitions, we assume in the following
that $S_1\coprod \ldots\coprod S_n=T_1\coprod \ldots\coprod T_k$.  The
proof of the second identity will be given in several steps.\par

Let us first introduce a useful notation. If $A$ is a twisted
algebra with product $m$ (respectively, a twisted coalgebra with
coproduct $\d$), $A^{\ot 2}$ and, more generally, $A^{\ot k}$ is
naturally provided with the structure of an algebra (respectively, of
a coalgebra). For example, the product on $A^{\ot k}$ (denoted by
$m_k$) is defined by:
\begin{eqnarray*}
  A^{\ot k}\ot A^{\ot k}
  & = & 
  (A_1\ot \cdots\ot A_k)\ot (A_1'\ot \cdots\ot A_k')\\[2mm]
  & \cong &
  (A_1\ot A_1')\ot \cdots\ot (A_k\ot A_k')\\[2mm]
  & \mapright{m^{\ot k}} &
  A^{\ot k},
\end{eqnarray*}
where we have written $A_i$ and $A_i'$ for copies of $A$. The
corresponding coproduct on $A^{\ot k}$ is denoted by $\d_k$.

\begin{lem} 
  Let $B$ be a twisted bialgebra. Then, for all $l$ and all $j$, we
  have:
  $$
  \d^{[l]}\circ m^{[j]}=m_l^{[j]}\circ (\d^{[l]})^{\ot j}.
  $$
\end{lem}

The proof may, for example, be given by induction on $l$ and is
essentially the same as the proof of the corresponding identity for
classical bialgebras (see \cite{pt} or \cite[Lem. II,8.]{p}).\par\ 
\par

The preceding lemma implies:
\begin{eqnarray*}
  X
  & := &
  1_{(S_1,\ldots,S_n)}\circ 1_{(T_1,\ldots,T_k)} \\[2mm]
  & = &
  m^{[n]}\circ 1^\ot_{(S_1,\ldots,S_n)}\circ\d^{[n]}
  \circ 
  m^{[k]}\circ 1^\ot_{(T_1,\ldots,T_k)}\circ\d^{[k]} \\[2mm]
  & = &
  m^{[n]}\circ 1^\ot_{(S_1,\ldots,S_n)}\circ m_n^{[k]}
  \circ
  (\d^{[n]})^{\ot k}\circ 1^\ot_{(T_1,\ldots,T_k)}\circ\d^{[k]}.
\end{eqnarray*}\par\ \par

\begin{lem} \label{copr-rule}
  We have, for $T\in\cal P$:
  $$
  \d\circ 1_T=\sum\limits_{U\coprod V=T}(1_U\ot 1_V)\circ\d
  $$ 
  and, more generally:
  $$
  \d^{[n]}\circ 1_T
  =
  \sum\limits_{U_1\coprod \ldots\coprod U_n=T}
   1^\ot_{(U_1,\ldots,U_n)}\circ\d^{[n]}.
  $$
\end{lem}

This follows from the definition of $\d$: indeed, $\d$ is a $\cal
P$-graded map, thus maps the degree $T$ component of $B$ to the degree
$T$ component of $B\ot B$, that is to
$$
\bigoplus\limits_{U\coprod V=T}B(U)\ot B(V).
$$
The same argument applies to $m$ and $m^{[k]}$. In particular:

\begin{lem} 
  We have:
  $$
      1^\ot_{(S_1,\ldots,S_n)}\circ m_n^{[k]}
     = 
    m_n^{[k]}
    \circ
    \sum\limits_{V_i^1\coprod \ldots\coprod V_i^k=S_i}
    1^\ot_{(V_1^1,\ldots,V_n^1,\;\ldots\;,V_1^k,\ldots,V_n^k)}.
$$
\end{lem}
 
Combined with Lemma~\ref{copr-rule}, this yields:
\begin{eqnarray*}
\lefteqn{%
X = m^{[n]}\circ m_n^{[k]}}\\[1mm]
& \circ & 
\Big(\sum\limits_{V_i^1\coprod \ldots\coprod V_i^k=S_i}
    1^\ot_{(V_1^1,\ldots,V_n^1,\;\ldots\;,V_1^k,\ldots,V_n^k)}\Big)\\[2mm]
& \circ &
\Big(\sum\limits_{U_1^i\coprod \ldots\coprod U_n^i=T_i}
    1^\ot_{(U_1^1,\ldots,U_n^1,\;\ldots\;,U_1^k,\ldots,U_n^k)}\Big)
\circ
(\d^{[n]})^{\ot k}\circ\d^{[k]}.  
\end{eqnarray*}
Besides, since the coproduct is coassociative, we have:
$(\d^{[n]})^{\ot k}\circ\d^{[k]}=\d^{[nk]}$ and, since the $1_S$ are
characteristic maps, we also have:
$$
    1^\ot_{(V_1^1,\ldots,V_n^1,\;\ldots\;,V_1^k,\ldots,V_n^k)}
\circ 
    1^\ot_{(U_1^1,\ldots,U_n^1,\;\ldots\;,U_1^k,\ldots,U_n^k)}
=0
$$
unless $V_i^j=U_i^j$ for all $i,j$. These identities are true if
and only if $V_i^j=U_i^j=S_i\cap T_j$. Therefore, the previous
expression of $X$ simplifies to:
$$
X
=
m^{[n]}\circ m_n^{[k]}
\circ 
1^\ot_{({S_1\cap T_1},\ldots,{S_n\cap T_1},\;
        \ldots\;, 
        {S_1\cap T_k},\ldots,{S_n\cap T_k})}
\circ 
\d^{[nk]}.
$$
The last point to notice is that, since the coproduct is
cocommutative, we can change the order of the tensor products of
characteristic maps freely in this expression, provided we modify
accordingly the multiplication map on the left. For example, because
of cocommutativity, the following identity holds:
$$
m_2^{[2]}\circ (1_U\ot 1_V\ot 1_W\ot 1_Y)\ot \d^{[4]}
=
(m^{[2]}\ot m^{[2]})\circ (1_U\ot 1_W\ot 1_V\ot 1_Y)\circ \d^{[4]},
$$
as well as all the general identities that can be constructed on
this pattern. In conclusion, we have
$$
X
=
m^{[n]}\circ (m^{[k]})^{\ot n}
\circ 
1^\ot_{({S_1\cap T_1},\ldots,{S_1\cap T_k},\;
        \ldots\;, 
        {S_n\cap T_1},\ldots,{S_n\cap T_k})}
\circ 
\d^{[nk]}.
$$
The proof of the theorem is complete upon noting that associativity
of $m$ implies $m^{[n]}\circ (m^{[k]})^{\ot n}=m^{[nk]}$. \par\ \par

The same computation, or a duality argument, show that, in case 
commutativity of $B$ was assumed, the second identity satisfied by the
elements of the twisted descent algebra would read:
$$
1_{( T_1,\ldots, T_k)} \circ 1_{( S_1,\ldots, S_n)} 
=
1_{(S_1\cap T_1,\ldots,S_1\cap T_k,\;\ldots\;,S_n\cap T_1,\ldots,S_n\cap T_k)}.
$$\par\ \par In other words, the structure of twisted descent algebras
of cocommutative twisted bialgebras is dual to the one of twisted
descent algebras of commutative twisted bialgebras.\par We conclude
this section with a remarkable identity for the linear generators of
${\cal T}_B$, relating the two associative products on twisted descent
algebras.

\begin{prop} \label{remarkable}
  Let $f,g,h,k$ be convolution products of characteristic maps
  in ${\cal T}_B$ such that $(f\circ g)\ast (h\circ k)\neq 0$.
  Then we have:
  $$
  (f\circ g)\ast (h\circ k) = (f\ast h)\circ (g\ast k).
  $$
\end{prop} 

Indeed, suppose that $f=1_{(S_1,\ldots,S_n)}$,
$g=1_{(T_1,\ldots,T_m)}$, $h=1_{(U_1,\ldots,U_p)}$,
$k=1_{(V_1,\ldots,V_q)}$.  Then the assumption $(f\circ g)\ast (h\circ
k)\neq 0$ implies $S_1\coprod \ldots\coprod S_n=T_1\coprod
\ldots\coprod T_m $, $U_1\coprod \ldots\coprod U_p=V_1\coprod
\ldots\coprod V_q$, and furthermore $S_i\cap U_j=\emptyset$ and
$T_i\cap V_j=\emptyset$ for all $i$, $j$.  The claim thus follows from
Theorem~\ref{main}.

\section{Solomon-Tits and the free twisted descent algebra.}

Let $(\Tau,\ast)$ be the algebra defined by generators and relations
as follows. The generators are the symbols $1_{ S}$, $ S\in\cal P$,
subject to the relations:
$$
1_{S_1}\ast \cdots\ast 1_{S_n}=0
$$
if $S_i\cap S_j\neq\emptyset$ for some $i\neq j$, and:
$$
1_\emptyset\ast 1_S = 1_S\ast 1_\emptyset =1_S.
$$
Note that the twisted descent algebra $\T_B$ of any twisted
bialgebra $B$ is isomorphic to a quotient of $\T$, according to our
previous computations.

\begin{prop} 
  The algebra $\Tau$ is, up to a canonical isomorphism, the twisted
  descent algebra of the twisted bialgebra freely generated, as a
  twisted algebra, by the primitive elements $\alpha_{ S}$, $ S\in\cal
  P$.
\end{prop}

Let us write $\cal B$ for the free twisted algebra on the generators
$\alpha_S,\ S\in\cal P$. A bijection from $S$ to $T$ maps $\alpha_S$
to $\alpha_T$ when $\cal B$ is viewed as a ${\cal P}$-graded vector
space. Then, since $\cal B$ is free:
$$
{\cal B}
=
\bigoplus\limits_{n\in\NM\cup\{0\}}
 (\bigoplus\limits_{S\in{\cal P}} {\bf k}\langle\alpha_S\rangle)^{\ot n},
$$
where we write ${\bf k}\langle\alpha_S\rangle$ for the vector
space generated by $\alpha_S$. It follows that ${\cal B}(S)$ is
spanned freely as a vector space by $\alpha_S$ and by the products
$\alpha_{S_1}\cdots\alpha_{S_k}$, where $S_1\coprod \ldots\coprod
S_k=S$, $S_i\not=\emptyset$, and where we write
$\alpha_{S_1}\cdots\alpha_{S_k}$ for the image in ${\cal B}(S)$ of
$\alpha_{S_1}\ot \cdots\ot\alpha_{S_k}\in {\cal B}(S_1)\ot
\cdots\ot{\cal B}(S_k)$.\par As a consequence, there is a unique
coproduct $\delta$ on ${\cal B}(S)$, defined by requiring that the
$\alpha_S$ are primitive elements:
$$
\delta (\alpha_{S_1}\cdots\alpha_{S_k})
=
\sum\limits_{\{i_1,\ldots,i_l\}\coprod\{j_1,\ldots,j_{k-l}\}=[k]}
  \alpha_{S_{i_1}}\cdots\alpha_{S_{i_l}}
  \ot 
  \alpha_{S_{j_1}}\cdots\alpha_{S_{j_{k-l}}},
$$
where $[k]=\{1,\ldots,k\}$.\par
Since the twisted descent algebra $\Tau_{\cal B}$ of $\cal B$ is a
quotient algebra of $\Tau$, the proposition is equivalent to the
following: the convolution product $1_{S_1}\ast \cdots\ast 1_{S_k}$,
$S_i\cap S_j= \emptyset\neq S_i$ are linearly independent in
$\Tau_{\cal B}$. \par 
To see this, note first the following consequence of the definition of
the convolution product in $\Tau_{\cal B}$ and of the coproduct in
$\cal B$. If $l\ge k$, then
$$
1_{(U_1,\ldots,U_l)}(\alpha_{S_1}\cdots\alpha_{S_k})
=
0
$$
unless $l=k$ and $(U_1,\ldots,U_l)$ may be obtained by reordering
$(S_1,\ldots,S_k)$; and in this case,
$$
1_{(U_1,\ldots,U_k)}(\alpha_{S_1}\cdots\alpha_{S_k})
=
\alpha_{U_1}\cdots\alpha_{U_k}.
$$

Let us assume now that 
$ f 
= \sum\limits_{k=1}^n
(\sum\limits_{i=1}^{m_k}\lambda_{i,k}1_{(S_1^{i,k},\ldots,S_k^{i,k})})=0
$
and that there exists an index $k$ such that $\lambda_{i,k}\neq 0$ for
some $i$. Choose $k$ minimal with this property. Then
$$
f(\alpha_{S_1}^{i,k}\cdots\alpha_{S_k}^{i,k})
=
\lambda_{i,k}\alpha_{S_1^{i,k}}\cdots\alpha_{S_k^{i,k}}+\Gamma ,
$$
where $\Gamma$ is a linear combination of noncommutative monomials
$\neq \alpha_{S_1^{i,k}}\cdots\alpha_{S_k^{i,k}}$ in the elements
$\alpha_{S_j^{i,k}}$ ($j\in[k]$), thus linearly independent of
$\alpha_{S_1^{i,k}}\cdots\alpha_{S_k^{i,k}}$. It follows that
$\lambda_{i,k}=0$, a contradiction. The proof of the proposition is
complete.\par\ \par

We will refer to the algebra $\Tau$ simply as ``the'' twisted descent
algebra (or, sometimes, as we did in the introduction, as the free
twisted descent algebra, when we want to emphasize that all descent
algebras of cocommutative twisted bialgebras are quotients of $\Tau$).
Note that, alternatively (and better suited for certain applications),
the twisted algebra $\Tau$ may be described as the free associative
algebra in the tensor category of connected $\cal P$-graded vector
spaces, with one generator $1_S$ in each degree. In general, we will
call a $\cal P$-graded algebra any algebra in the tensor category of
connected $\cal P$-graded vector spaces, so that $\Tau$ is a free
$\cal P$-graded algebra. That is, if $K=\bigoplus\limits_{S\in {\cal
    P}, S\neq\emptyset} {\bf k}$, we have:
$$
\Tau 
=
\bigoplus\limits_{n\in \NM\cup\{0\}} K^{\ot n}
=
\bigoplus\limits_{(S_1,\ldots,S_n)\in{\cal P}^n,S_i\not=\emptyset} {\bf k}.
$$
The component of degree $S$ ($S\in\cal P$) of $\Tau$ is:
$$
\Tau_S
=
\bigoplus\limits_{\ S_1\coprod\ldots\coprod S_n=S,\ S_i\not=\emptyset} {\bf k},
$$
and $\Tau_\emptyset ={\bf k}^{\ot 0}=\bf k$.  
The $\cal P$-graded product is given on the indices by:
$$
(( S_1,\ldots, S_n),(T_1, \ldots, T_l))
\longmapsto 
( S_1,\ldots, S_n, T_1, \ldots, T_l).
$$
%% In what follows, we write $1_{( S_1,\ldots, S_n)}$ for $1_{S_1}\ast
%% \cdots\ast 1_{S_n}\in\T$ again, for all $S_i\in{\cal P}$.  
Recall that, by
Theorem~\ref{main}, there is a second associative product $\circ$ on
$\Tau$, the composition product, which is homogeneous with respect to
the $\cal P$-grading. That is,
$$
1_{( S_1,\ldots, S_n)} \circ 1_{( T_1,\ldots, T_k)}=0
$$
if $S_1\coprod \ldots\coprod S_n\not= T_1\coprod \ldots\coprod
T_k$. In particular, the algebra $(\Tau ,\circ )$ splits as a product
of its graded components.  In other words, the graded component
$\Tau_S$ of $\Tau$ of degree $S\in\cal P$ is a two-sided ideal of
$\Tau$ with respect to $\circ$. The algebraic structure of this ideal
is studied in detail in \cite{sch2}.  Note that, if $S,T\in{\cal P}$
such that $|S|=|T|$, then any bijection from $S$ to $T$ induces an
isomorphism of algebras from $(\Tau_S,\circ)$ to $(\Tau_T,\circ)$.

\begin{prop} 
  The algebra $\Tau_{[n]}$ is isomorphic to the Solomon-Tits
  algebra of the symmetric group on $[n]$.
\end{prop}

Recall briefly the definitions in Tits' appendix to Solomon's original
paper \cite{s}. Let $(W,S)$ be a Coxeter system (that is, $W$ is a
finite Coxeter group and $S$ a simple system of generators, see e.g.
\cite{hum}). For $K\subset S$, the subgroup of $W$ generated by $K$ is
written $W_K$. The Coxeter complex $\Sigma$ associated to $(W,S)$ is a
simplicial complex. Its simplices are in bijection with the cosets
$W_K\cdot w$ with $K\subset S$ and $w\in W$. The fundamental
geometrical idea, due to Tits \cite{ti2}, is that there is a product
on $\Sigma$, the associativity of which is proven in \cite{ti}. The
product of $A$ and $B$ in $\Sigma$ is, by definition, ``the greatest
common face of the first terms of all galleries of minimum length
whose first and last term dominate, respectively, $A$ and $B$''
(recall that a gallery is a sequence of adjacent maximal simplices in
the Coxeter complex). \par We call the algebra $\Sigma$, for this
geometrical product, the Solomon-Tits algebra, since it was introduced
by Tits in order to explain geometrically the meaning of Solomon's
constructions in \cite{s}, and to relate them to Tits' previous results on
buildings. The connection to the algebra $\Tau_{[n]}$ is
not obvious, but becomes so if one chooses the right parametrization
of the simplices of the Coxeter complex of $S_n$. We learned this
connection in Brown \cite{br}, who also studies various
generalizations and applications of Tits' calculations on hyperplane
arrangements.\par The Coxeter complex can be viewed as a triangulation
of the sphere associated to the hyperplane arrangement corresponding
to $S_n$, see \cite{hum} for details. This arrangement is obtained
from the hyperplanes $x_i=x_j,\ i\not=j$ in $\RM^n$ ($S_n$ is not
essential relative to $\RM^n$, leaving a $1$-dimensional subspace
invariant, but this does not matter at this point). The maximal cells
of the corresponding triangulation of the sphere correspond to the
$n$-dimensional simplicial cones in $\RM^n$ obtained from the
arrangement. That is, since the half-spaces associated to the
hyperplanes are given by the condition $x_i\geq x_j$, there are $n!$
such cones, described by the following set of inequalities
in the canonical coordinate system:
$$
x_{\s (1)}\geq \cdots\geq x_{\s (n)},
$$
where $\s\in S_n$. The faces of the cones are obtained by turning
an arbitrary number of inequalities into equalities. For example,
$x_3\geq x_1=x_2\geq x_4$ is a face of both $x_3\geq x_2\geq x_1\geq
x_4$ and $x_3\geq x_1\geq x_2\geq x_4$. It follows in particular from
this description of $\Sigma$ that its elements are parametrized by the
ordered partitions of $[n]$; for instance, the partition $(\{3\},
\{1,2\}, \{4\})$ is associated to the face $x_3\geq x_1=x_2\geq x_4$.
Using this dictionary, it is an easy exercise to express Tits' product
on $\Sigma$ by using the defining formulas for $\circ$, see \cite{br}.
The proposition follows.\par\ \par

As a final structural ingredient, in the free case, we have:

\begin{prop} 
  The algebra $\Tau$ carries a coassociative and cocommutative
  coproduct, defined by:
  $$
  \d (1_{( S_1,\ldots, S_k)})
  :=
  \sum\limits_{T_i\coprod U_i=S_i} 
  1_{(T_1,\ldots,T_k)}\ot 1_{(U_1,\ldots,U_k)}.
  $$
\end{prop}

Empty sets are allowed in this summation formula. The coassociativity
of the coproduct is straightforward. The same formula defines a
cocommutative $\cal P$-graded coalgebra structure on $\Tau$.\par

\begin{thm} \label{bialgebra}
  The product $\ast$ and the coproduct $\d$ turn $\Tau$ into a $\cal
  P$-graded bialgebra. The product $\circ$ and the coproduct $\d$ turn
  $\Tau$ into a bialgebra.
\end{thm}

Once again, this property parallels the fundamental properties of the
classical descent algebra, see e.g. \cite{pr2}. Notice, however, that
the $\cal P$-graded hypothesis is necessary, since $(\Tau,\ast,\d)$ is
not bialgebra. For example,
$$
\d (1_{\{1,2\}}\ast 1_{\{1,2\}})=\d (0)=0,
$$
while
\begin{eqnarray*}
  \lefteqn{%
    \d(1_{\{ 1,2\} })\ast_2 \d (1_{\{ 1,2\} })}\\[2mm]
  & = &
  \mbox{\hphantom{$\ast_2$}}
  (1_{\{ 1,2\} }\ot 1_\emptyset 
  +1_{\{ 1\} }  \ot 1_{\{ 2\} }
  +1_{\{ 2\} }  \ot 1_{\{ 1\} }
  +1_\emptyset  \ot 1_{\{ 1,2\} })\\[1mm]
  & & 
  \ast_2
  (1_{\{ 1,2\} }\ot 1_\emptyset +1_{\{ 1\} }\ot 1_{\{ 2\}
  }+1_{\{ 2\} }\ot 1_{\{ 1\} }+1_\emptyset\ot 1_{\{ 1,2\} })\\[2mm]
  & = &
  2\cdot 1_{\{ 1,2\} }       \ot 1_{\{ 1,2\} }
  +1_{(\{ 1\},\{2\} )} \ot 1_{(\{ 2\},\{1\} )}
  +1_{(\{ 2\},\{1\} )} \ot 1_{(\{ 1\},\{2\} )}.
\end{eqnarray*}\par
The first part of the theorem follows from the following identity,
that holds for all mutually disjoint sets 
$ S_1,\ldots,S_n, V_1,\ldots  ,V_k\in{\cal P}$:
\begin{eqnarray*}
\lefteqn{%
  \d(1_{( S_1, \ldots,S_n)}\ast 1_{ ( V_1,\ldots, V_k)})}\\[3mm]
& = &
\d(1_{(S_1,\ldots,S_n,V_1,\ldots,V_k)}) \\[2mm]
& = &
\sum\limits_{ T_i\coprod U_i= S_i, W_j\coprod Z_j= V_j}
 1_{(T_1,\ldots,T_n,W_1,\ldots,W_k)}
\ot 
1_{(U_1,\ldots,U_n,Z_1,\ldots,Z_k)} \\[2mm]
& = &
\sum\limits_{ T_i\coprod U_i= S_i, W_j\coprod Z_j= V_j}
 [1_{(T_1,\ldots,T_n)}\ot 1_{(U_1,\ldots,U_n)}]
\ast_2
 [1_{(W_1,\ldots,W_k)}\ot 1_{(Z_1,\ldots,Z_k)}]\\[3mm]
& = &
\d (1_{(S_1,\ldots,S_n)})\ast_2 \d (1_{(V_1,\ldots,V_k)}).  
\end{eqnarray*}
Let us assume now that $S_1\coprod \ldots\coprod
S_n=V_1\coprod \ldots\coprod V_k$. Then, on the one hand, we have:
\begin{eqnarray*}
\lefteqn{%
A
 := 
\d (1_{( S_1, \ldots,S_n)}\circ 1_{ ( V_1,\ldots, V_k)})}\\[3mm]
&  = &
\d (1_{( {S_1\cap V_1} ,\ldots, {S_1\cap V_k},\;\ldots\;, {S_n\cap
    V_1},\ldots, {S_n\cap V_k})})\\[2mm]
&  = &
\sum\limits_{U_{i,j}\coprod T_{i,j}=S_i\cap V_j}
1_{(U_{1,1},\ldots,U_{1,k},\;\ldots\;,U_{n,1},\ldots,U_{n,k})}
\ot 1_{(T_{1,1},\ldots,T_{1,k},\;\ldots\;,T_{n,1},\ldots,T_{n,k})}.
\end{eqnarray*}
On the other hand, we have:
\begin{eqnarray*}
\lefteqn{%
B
 := 
\d (1_{( S_1, \ldots,S_n)})\circ_2 \d (1_{ ( V_1,\ldots, V_k)})}\\[2mm]
&  = &
[\sum\limits_{U_{i}\coprod T_{i}=S_i}
1_{(U_{1},\ldots,U_{n})}\ot 1_{(T_{1},\ldots,T_{n})}]
 \circ_2 
[\sum\limits_{U^{j}\coprod T^{j}=V_j}
1_{(U^{1},\ldots,U^{k})}\ot 1_{(T^{1},\ldots,T^k)}]\\[2mm]
&  = &
\sum\limits_{U_{i}\coprod T_{i}=S_i,\ U^{j}\coprod T^{j}=V_j}
 [1_{(U_{1},\ldots,U_{n})}\circ (1_{(U^{1},\ldots,U^{k})}]
\ot 
[1_{(T_{1},\ldots,T_{n})}\circ 1_{(T^{1},\ldots,T^k)}].  
\end{eqnarray*}
Each term in the last sum is $0$ unless $U_1\coprod \ldots\coprod
U_n=U^1\coprod \ldots\coprod U^k$ and $T_1\coprod \ldots\coprod
T_n=T^1\coprod \ldots\coprod T^k$. In particular, in this case,
it follows that $U_i\cap U^j=U_i\cap V_j$. For otherwise, the
elements of $U_i\cap V_j$ not in $U_i\cap U^j$ would lie in $T^j$ and
would not belong to any $T_l$, and the term would be $0$, a
contradiction. In conclusion, setting $U_{i,j}:=U_i\cap U^j$ and
$T_{i,j}:=T_i\cap T^j$, yields the requested identity, $A=B$, and
completes the proof of Theorem~\ref{bialgebra}.
\par The existence of the coproduct $\delta$ allows 
to derive the following variant of Proposition~\ref{remarkable}
for the algebra $\Tau$, now relating both products and the coproduct
on the free twisted descent algebra.

\begin{cor} \label{reciproci}
  Let $f,g,h\in{\cal T}$, then we have:
  $$
  (f\ast g)\circ h = m((f\ot g)\circ_2 \d(h)),
  $$
  where $m:{\cal T}\ot{\cal T}\to {\cal T}$ is the convolution product.
\end{cor} 
In more illustrative terms, using Sweedler's notation
$\delta(h)=\sum h^{(1)}\ot h^{(2)}$ for the coproduct, we have
$$
(f\ast g)\circ h = \sum (f\circ h^{(1)})\ast (g\circ h^{(2)}).
$$
For the proof, it suffices to consider $f=1_{(S_1,\ldots,S_n)}$,
$g=1_{(T_1,\ldots,T_k)}$ and $h=1_{(U_1,\ldots,U_p)}$, by linearity.
Then, the right hand side
\begin{eqnarray*}
  \lefteqn{%
    m((f\ot g)\circ_2 \d(h))}\\[2mm]
  & = &
  \sum\limits_{X_i\coprod Y_i=U_i}
   m(
   (1_{(S_1,\ldots,S_n)}\ot 1_{(T_1,\ldots,T_k)})
   \circ_2
   (1_{(X_1,\ldots,X_p)}\ot 1_{(Y_1,\ldots,Y_p)}))\\[2mm]
  & = &
  \sum\limits_{X_i\coprod Y_i=U_i}
   (1_{(S_1,\ldots,S_n)}\circ 1_{(X_1,\ldots,X_p)})
   \ast
   (1_{(T_1,\ldots,T_k)}\circ 1_{(Y_1,\ldots,Y_p)})
\end{eqnarray*}
does not vanish if and only if $S_i\cap T_j=\emptyset$ for all $i$,
$j$ and $ S\coprod T=U$, where $S=S_1\coprod \ldots\coprod S_n$,
$T=T_1\coprod \ldots\coprod T_k$ and $U=U_1\coprod \ldots\coprod U_p$.
The same observation is true for the left hand side.  And in this
case, there is a unique summand $\neq 0$ in the above sum, indexed by
$X_i=S\cap U_i$ and $Y_i=T\cap U_i$, that is:
\begin{eqnarray*}
  \lefteqn{%
    m((f\ot g)\circ_2 \d(h))}\\[2mm]
  & = &
   (1_{(S_1,\ldots,S_n)}\circ 1_{(S\cap U_1,\ldots,S\cap U_p)})
   \ast
   (1_{(T_1,\ldots,T_k)}\circ 1_{(T\cap U_1,\ldots,T\cap U_p)})\\[2mm]
   & = &
    1_{(S_1\cap U_1,\ldots,S_1\cap U_p,
        \;\ldots\;,
        S_n\cap U_1,\ldots,S_n\cap U_p,
        T_1\cap U_1,\ldots,T_1\cap U_p,
        \;\ldots\;,
        T_k\cap U_1,\ldots,T_k\cap U_p)}\\[2mm]
   & = &
(f\ast g)\circ h
\end{eqnarray*}
as asserted.\par\ \par 

Recall that the (ordinary) descent algebra of a twisted bialgebra $B$
(as introduced in \cite{pr3}) may be identified with the convolution
algebra generated by the formal series $1_n:=\sum\limits_{S\in{\cal
    P},\ |S|=n}1_S$, $n\in\NM$.  Motivated by this remark, in the free
case, we write $\cal D$ for the convolution subalgebra of $\hat \Tau
:=\prod\limits_{S\in\cal P}\Tau_S$ generated by the elements
$1_n\in\hat{\Tau}$, $n\in\NM$.

\begin{thm} \label{desc-alg}
  The convolution algebra ${\cal D}$ is also a subalgebra of
  $\hat\Tau$ with respect to the composition product. More precisely,
  each graded component
  $$
  {\cal D}_n:={\cal D}\cap \prod\limits_{|S|=n}\Tau_S
  $$
  is closed under the composition $\circ$ of endomorphisms of ${\cal B}$.
  
  Moreover, $({\cal D}_n,\circ )$ is isomorphic to Solomon's classical
  descent algebra of the symmetric group on $[n]$.
\end{thm}

Notice that, as a particular consequence of the theorem,
Corollary~\ref{reciproci} is a generalization of the crucial
reciprocity law for the descent algebra derived in
\cite[Proposition~5.2]{g}.

Concerning the proof, we observe that ${\cal D}_n\circ{\cal D}_m=0$
whenever $n\neq m$, and that the elements
$1_{n_1,\ldots,n_k}:=1_{n_1}\ast \cdots \ast 1_{n_k}$ with
$n_1+\cdots+n_k=n$ constitute a linear basis of ${\cal D}_n$.
It remains to prove Solomon's fundamental multiplication rule \cite[Theorem~1]{s}
for the members of this basis, that is:
$$
1_{n_1,\ldots,n_k}\circ 1_{m_1,\ldots,m_l}
=
\sum 1_{a_1^1,\ldots,a_1^l,\;\ldots\;,a_k^1,\ldots,a_k^l},
$$
where the sum on the right is taken over all matrices $(a_i^j)$ of
nonnegative integers such that $\sum\limits_{j=1}^l a_i^j=n_i$ and
$\sum\limits_{i=1}^k a_i^j=m_j$ (see \cite{r,p}).

By definition, we have
$$
1_{n_1,\ldots,n_k}
=
\sum\limits_{|S_i|=n_i,\ S_i\cap S_j=\emptyset} 1_{(S_1,\ldots,S_k)}
$$
and therefore:
\begin{eqnarray*}
  \lefteqn{%
    1_{n_1,\ldots,n_k}\circ 1_{m_1,\ldots,m_l}}\\[2mm]
  & = &
  \sum\limits_{|S_i|=n_i,\ S_i\cap S_j=\emptyset}\quad
   \sum\limits_{|T_u|=m_u,\ T_u\cap T_v=\emptyset}
    1_{(S_1,\ldots,S_k)}\circ 1_{(T_1,\ldots,T_l)}\\[2mm]
  & = &
  \sum\limits_{%
    |S_i|=n_i,|T_u|=m_u,
    S_1\coprod \ldots\coprod S_k=T_1\coprod \ldots\coprod T_l}
  1_{%
   (S_1\cap T_1,\ldots,S_1\cap T_l,\;\ldots\;,S_k\cap T_1,\ldots,S_k\cap T_l)}.
\end{eqnarray*}
Writing $U_i^j$ for $S_i\cap T_j$, the last term reads:
\begin{eqnarray*}
  & = &
  \sum\limits_{\sum\limits_{j=1}^l a_i^j=n_i,\sum\limits_{i=1}^ka_i^j=m_j}
  \sum\limits_{|U_i^j|=a_i^j}
  1_{(U_1^1,\ldots,U_1^l,\ldots,U_k^1,\ldots,U_k^l)},
  \hspace*{25ex}
\end{eqnarray*}
and the theorem is proven.

\section{Applications to computations with shuffles}

In this last section, we derive some important identities in the
twisted descent algebra ${\cal T}$. They hold in the twisted descent
algebra of any twisted bialgebra $B$, since $\Tau_B$ is a quotient of
$\Tau$. When these identities follow immediately from
Theorem~\ref{main}, the proofs are omitted. \par

\begin{lem} 
  Let $S=\{s_1,\ldots,s_n\}\in\cal P$. Then, for all 
$S_1,\ldots,S_k$ in $\cal P$ such that $S_1\coprod
  \ldots\coprod S_k=S$, we have:
  $$
  1_{(\{s_1\},\ldots,\{s_n\})}
  \circ
  1_{(S_1,\ldots,S_k)}
  =
  1_{(\{s_1\},\ldots,\{s_n\})}.
  $$
\end{lem}

The next lemma shows that the behavior of products in $(\Tau,\circ)$ encodes
the unshuffling of permutations, an important device in algebraic
combinatorics and its applications (card shuffling, models of
databases structure in computer science, and so on).

\begin{lem} 
  Let $S$ and $S_1,\ldots,S_k$ be as above. Then we have
  $$
  1_{(S_1,\ldots,S_k)}\circ 1_{(\{s_1\},\ldots,\{s_n\})}=1_{\b (S)},
  $$
  where $\b (S)$ is the $(S_1,\ldots,S_k)$-relative unshuffling of
  $(\{s_1\},\ldots,\{s_n\})$.
\end{lem} 

That is, $\b (S)$ is the sequence of one-element subsets of $S$
obtained by selecting successively in the sequence 
$(\{s_1\},\ldots,\{s_n\})$ the
elements of $S_1$, $S_2$, \ldots, $S_k$. For example, we have:
$$
1_{(\{1,3,5\},\{2,4\})}\circ 1_{(\{3\},\{4\},\{5\},\{2\},\{1\})}
=
1_{(\{3\},\{5\},\{1\},\{4\},\{2\})}.
$$
>From now on, we write $1_\s$ for $1_{(\{\s(1)\},\ldots,\{\s (n)\})}$
whenever $\sigma\in S_n$.

\begin{cor} \label{shuffle}
  Let $(S_1,\ldots,S_k)$ be an increasing partition of $[n]$. That is,
  we assume that $s<s'$ for all $s\in S_i,\ s'\in S_j$ whenever $i<j$.
  Then,
  $$
  1_{(S_1,\ldots,S_k)}\circ 1_{\s}=1_{(\{1\},\ldots,\{n\})}
  $$
  if and only if $\s$ is a shuffle of $S_1,\ldots,S_k$.
\end{cor}

Recall that a finite sequence $C=(n_1,\ldots,n_k)$ of positive
integers with sum $n$ is a composition of $n$.  Recall furthermore
that the descent set $D(\s)$ of a permutation $\s\in S_n$ is the
subset of $[n-1]$ consisting of all $i$ such that $\s(i)>\s(i+1)$.  We
denote by $D_C$ the sum of all permutations in $S_n$ with descent set
contained in $\{n_1,n_1+n_2,\ldots,n_1+\cdots+n_{k-1}\}$.\par
The linear span of these elements in the integral group ring $\ZM
[S_n]$ is the classical descent algebra of $S_n$, in Solomon's
original setting.  It is closed under composition of permutations, due
to Solomon's remarkable discovery. The elements $D_C$ form a linear
basis of Solomon's algebra, and the corresponding structure
coefficients (or, equivalently, those of the dual algebra) can be
computed explicitly in terms of double coset representatives of Young
subgroups of $S_n$, or using Tits' algebraic approach to the
geometrical properties of Coxeter complexes. An equivalent
combinatorial description already occured in the proof of
Theorem~\ref{desc-alg}.\par\ \par If $(S_1,\ldots,S_k)$ is a partition
of $[n]$, we say that the composition $C=(|S_1|,|S_2|,\ldots,|S_k|)$
of $n$ is the type of $(S_1,\ldots, S_k)$.

\begin{prop} 
  Let $(S_1,\ldots,S_k)$ be an increasing partition of $[n]$ of type~$C$. 
  Then the sum of all permutations $\s \in S_n$
  such that
  $$
  1_{(S_1,\ldots,S_k)}\circ 1_{\s}=1_{(\{1\},\ldots,\{n\})}
  $$
  is identical to $D_C^\ast$,
  where $\ast$ is the linear endomorphism of $\ZM[S_n]$ defined by $\s^\ast
  :=\s^{-1}$ for all $\s\in S_n$.
\end{prop} 

This is another way of stating Corollary~\ref{shuffle}, since the
inverse of a shuffle of $S_1,\ldots,S_k$ is (when viewed as a word) a
product $u_1\ldots u_k$ of increasing words $u_i$ of length $|S_i|$.
In other terms, the mapping $\s\mapsto \s^{-1}$ yields a bijection
from the $(S_1,\ldots,S_k)$-shuffles onto the set of all $\s\in S_n$
such that $D(\s)\subseteq
\{|S_1|,|S_1|+|S_2|,\ldots,|S_1|+\cdots+|S_{k-1}|\}$, as asserted.\par
\ \par

Recall that $\Tau_n:=\Tau_{[n]}$ (the graded component of $\Tau$ of
degree $[n]$) is an algebra for the product $\circ$, with linear basis
consisting of the elements $1_{(S_1,\ldots,S_k)}$ such that
$S_1\coprod \ldots\coprod S_k=[n]$ and $S_i\neq\emptyset$.  There is a
right action of $S_n$ on this basis, defined by:
$$
1_{(S_1,\ldots,S_k)}\cdot\s
:=
1_{(\s^{-1} (S_1),\ldots,\s^{-1} (S_k))},
$$
which turns $\Tau_n$ into a permutation module for $S_n$, by
linearity.

\begin{prop} 
  The product $\circ$ on $\Tau_n$ is equivariant with respect to the
  $S_n$-action. 
\end{prop} 

The proposition follows from the definition of $\circ$ since, for any
subsets $S$ and $T$ of $[n]$, we have:
$$
\s^{-1} (T)\cap \s^{-1} (S)=\s^{-1} (T\cap S).
$$
Note that $1_{(S_1,\ldots,S_k)}$ and $1_{(T_1,\ldots,T_l)}$ belong
to the same $S_n$-orbit if and only if $(S_1,\ldots,S_k)$ and
$(T_1,\ldots,T_l)$ are of the same type.
In particular, the fixed space ${\cal F}_n$ of $S_n$ in $\T_n$ is linearly
generated by the orbit sums
$$
O_{n_1,\ldots,n_k}
=
\sum_{|S_i|=n_i} 1_{(S_1,\ldots,S_k)},
$$
where $n_1+\cdots+n_k=n$.  Furthermore, $S_n$-equivariance of the
composition product on $\T_n$ implies that ${\cal F}_n$ is a
subalgebra of $(\T_n,\circ)$. The algebra ${\cal F}_n$ is also
isomorphic to Solomon's classical descent algebra. This surprising
observation is due to Bidigare (\cite{bi}, see also \cite{br,sch}).
The link to our considerations in the previous section is given by the
following simple fact. Each linear generator $1_{n_1,\ldots,n_p}$ of
${\cal D}_n$ is the sum of its graded components
$$
1_{n_1,\ldots,n_k}(S)
=
\sum\limits_{|S_i|=n_i,\ S_1\coprod\ldots\coprod S_k=S} 
     1_{(S_1,\ldots,S_k)} 
\qquad(|S|=n),
$$
and we have:

\begin{prop}
  The truncation mapping 
  $$
  1_{n_1,\ldots,n_k}
  = 
  \sum_{|S|=n} 1_{n_1,\ldots,n_k}(S)
  \longmapsto 
  1_{n_1,\ldots,n_k}([n])
  =
  O_{n_1,\ldots,n_k}
  $$
  is an isomorphism of algebras from $({\cal D}_n,\circ)$ onto
  $({\cal F}_n,\circ)$.
\end{prop}

This is immediate from the multiplication rule
$1_{n_1,\ldots,n_k}(S)\circ 1_{m_1,\ldots,m_l}(T)=0$ for all $S\neq
T$.\par\ \par
To conclude, observe that the stabilizer of $1_{(S_1,\ldots,S_k)}$ in
$S_n$ is a parabolic subgroup (i.e. it is conjugated to a Young
subgroup).  If, in particular, $(S_1,\ldots,S_k)$ is an increasing
partition of $[n]$ and $n_i=|S_i|$, then the stabilizer of
$1_{(S_1,\ldots,S_k)}$ is the Young subgroup $S_{n_1}\times
\cdots\times S_{n_k}\subset S_n$.  In other words, as an $S_n$-module,
$\T_n$ is the direct sum of the Young modules arising from
the action of $S_n$ on the left cosets of $S_{n_1}\times\cdots\times S_{n_k}$
in $S_n$, $n_1+\cdots+n_k=n$. These cosets are well-known to
be parametrized by $(S_1,\ldots,S_k)$-shuffles, that is, by the
summands of $D_C^*$.  More precisely, any permutation in $S_n$ can be
written uniquely as the composition of an element in the corresponding
Young subgroup with a $(S_1,\ldots,S_k)$-shuffle.\par 
This property can be recovered easily using the twisted descent
formalism: first, for any $\s\in S_n$, the corresponding element of
the Young subgroup is the permutation $\b$ defined by:
$$
1_{(S_1,\ldots,S_k)}\circ 1_{\s}=1_{\b }.
$$
The element $\b^{-1}\cdot \s$ is then a $(S_1,\ldots,S_k)$-shuffle.
Indeed, we have:
\begin{eqnarray*}
1_{(S_1,\ldots,S_k)}\circ 1_{\b^{-1}\cdot \s }
& = &
1_{(S_1,\ldots,S_k)}\circ (1_{\s }\cdot \b ) \\[2mm]
& = &
(1_{(S_1,\ldots,S_k)}\circ 1_{\s })\cdot \b,\hspace*{10ex}
\end{eqnarray*}
by equivariance of the $\circ$ product, and since $\b$ is an
element of the Young subgroup stabilizing $1_{(S_1,\ldots,S_k)}$, 
\begin{eqnarray*}
\mbox{\hphantom{$1_{(S_1,\ldots,S_k)}\circ 1_{\b^{-1}\cdot \s }$}}
& = &
1_{\b }\cdot \b =1_{\b^{-1}\cdot \b}=1_{(\{1\},\ldots,\{n\})}.
\end{eqnarray*}
Now Corollary~\ref{shuffle} implies the desired property.\par\ \par
It follows from these remarks, properties and identities, that there
are straight connections between the multiple algebraic structure on
$\T$ and the classical descent algebra ${\cal D}$.  In particular, as
already mentioned in the introduction, Solomon's computations that
involve fine computations with coset representatives of Young
subgroups \cite{s}, can be rephrased and handled in the language of
twisted descents in the same way that they were translated by Tits in
the language of Coxeter complexes and buildings (see \cite[pp.257]{s}
and \cite{ti}).\par


\begin{thebibliography}{99}
\bibitem{ba} M.G. Barratt. Twisted Lie algebras.
Geom. Appl. Homotopy Theory II, Proc. Conf., Evanston 1977, \it Lect. Notes Math. \rm 658, 9-15 (1978).
\bibitem{bi} P. Bidigare.
{H}yperplane {A}rrangement {F}ace {A}lgebras and their {A}ssociated
{M}arkov {C}hains.
Ph.d. thesis, {U}niversity of {M}ichigan (1997).
\bibitem{bid} P. Bidigare, P. Hanlon, D. Rockmore.
A combinatorial description of the spectrum for the Tsetlin library and its generalization to hyperplane arrangements. 
\it Duke Math. J. \rm 99, No.1, 135-174 (1999).
\bibitem{br} K. S. Brown,
Semigroups, rings, and Markov chains. 
\it J. Theor. Probab. \rm 13, No.3, 871-938 (2000).
\bibitem{ch} K.-T. Chen, \it Collected papers of K.-T. Chen. \rm Edited and with a preface by Philippe Tondeur, and an essay on Chen's life and work by Richard Hain and Tondeur. Contemporary Mathematicians.
 Birkh\"auser Boston, Inc., Boston, MA, 2001.
\bibitem{do} A. Dold.
\it Lectures on algebraic topology. \rm 2nd ed. 
Grundlehren der mathematischen Wissenschaften 200. Berlin Heidelberg New York: Springer-Verlag. XI, 1980.
\bibitem{g} I. M. Gelfand, D. Krob, A. Lascoux, B. Leclerc, V. Retakh,
J.-Y. Thibon. Noncommutative symmetric functions.\it Adv. Math. \rm 112, No.2,
218-348 (1995).
\bibitem{hum} J. Humphreys, \it Reflexion Groups and Coxeter Groups. \rm Cambridge University Press, 1990.
\bibitem{joy}A. Joyal.
Foncteurs analytiques et esp\`eces de structures.
Combinatoire \'enum\'erative, Proc. Colloq., Montr\'eal/Can. 1985, \it Lect. Notes Math. \rm 1234, 126-159 (1986).
\bibitem{k} D. Krob, B. Leclerc, J.-Y. Thibon.
Noncommutative symmetric functions. II: Transformations of alphabets. \it Int. J. Algebra Comput. \rm 7, No.2, 181-264 (1997).
\bibitem{mr} C. Malvenuto, C. Reutenauer.
Duality between quasi-symmetric functions and the Solomon descent algebra. \it J. Algebra \rm 177, No.3, 967-982 (1995).
\bibitem{p2} F. Patras, Adams operations, iterated integrals and free loop spaces cohomology. \it J. Pure Appl. Algebra. \rm 114, (1997), 273-286. 
\bibitem{pa} F. Patras. Adams operations, algebras up to homotopy and cyclic homology. \it Topology. \rm 39, (6), (2000), 1089-1101. 
\bibitem{p} F. Patras. L'alg\`ebre des descentes d'une big\`ebre
  gradu\'ee. \it J. Algebra. \rm 170, No.2, (1994), 547-566.
\bibitem{pt} F. Patras.
La d\'ecomposition en poids des alg\`ebres de Hopf. \it Ann. Inst. Fourier. \rm
43, No.4, 1067-1087 (1993).
\bibitem{pr2}F. Patras; C. Reutenauer. Lie representations and an algebra containing Solomon's. \it J. Alg. Comb.\rm , 16, (2002), 301-314. 
\bibitem{pr3}F. Patras; C. Reutenauer. On descent algebras and twisted bialgebras. \it Moscow Math. J. \rm To appear.
\bibitem{ra} G. Racinet, Doubles m\'elanges des polylogarithmes multiples aux racines de l'unit\'e. \it Publ. Math. Inst. Hautes \'Etudes Sci. \rm No. 95 (2002), 185--231.
\bibitem{r} C. Reutenauer. \it{Free Lie algebras}. \rm Oxford University Press, 1993.
\bibitem{sch} M. Schocker. The descent algebra of the symmetric group.
  To appear in the Proceedings of the Instructional Workshop of the
  10th International Conference on Representations of Algebras and
  Related Topics (ICRA~X), held at the Fields Institute, Toronto, July
  2002.
\bibitem{sch2} M. Schocker. Semigroups, Young modules, and the descent
  algebra of the symmetric group. Preprint.
\bibitem{s} L. Solomon.
A Mackey formula in the group algebra of a finite Coxeter group. \it J. Algebra
\rm 41, 255-268 (1976). 
\bibitem{st} C. R. Stover.
The equivalence of certain categories of twisted Lie and Hopf algebras over a commutative ring.
\it J. Pure Appl. Algebra \rm 86, No.3, 289-326 (1993).
\bibitem{sw} M. Sweedler. \it Hopf algebras. \rm Benjamin, New York, 1969.
\bibitem{ti} J. Tits. Two properties of Coxeter complexes. Appendix to: L. Solomon.
A Mackey formula in the group algebra of a finite Coxeter group. \it J. Algebra
\rm 41, 255-268 (1976).
\bibitem{ti2} J. Tits. \it Buildings of spherical type and finite $BN$-pairs\rm , Springer-Verlag, Berlin, 1974.
\end{thebibliography}
\end{document}